\documentclass[11pt,english,reqno]{amsart}

\RequirePackage[utf8,latin1]{inputenc}
\usepackage{tipa}
\usepackage{amstext,amsopn,amscd,amsthm}
\usepackage{amsfonts,amsmath,amssymb,amsbsy}
\usepackage{layout}
\usepackage{newlfont}
\usepackage{rotating}
\usepackage{verbatim}
\usepackage{lipsum}

\RequirePackage{graphicx}
\graphicspath{{Pictures/}}

\usepackage{caption}
\usepackage{subcaption}

\usepackage{picture}
\usepackage{wrapfig}
\usepackage[dvipsnames,svgnames,x11names]{xcolor}
\usepackage{pstcol,pst-fill,pst-grad}
\usepackage{lscape}

\usepackage{array}
\usepackage{diagbox}
\usepackage{colortbl}
\usepackage{hhline}
\usepackage{rotating,multirow,multicol}

\usepackage[bookmarks=false,colorlinks=true,linkcolor=DarkRed,urlcolor=DarkRed,citecolor=DarkRed,filecolor=DarkRed,linktocpage=true,pdfstartview={XYZ null null 0.702}]{hyperref}

\usepackage{titletoc}
\newcommand{\noi}{\noindent}\newcommand{\ppass}{\vspace{0.2cm}\noindent}

\title[Who invented von Koch's snowflake curve?]{Who invented von Koch's snowflake curve?}

\author[Yann Demichel]{\href{mailto:yann.demichel@parisnanterre.fr}{Yann Demichel\textsuperscript{\dag}}}

\address{\textsuperscript{\dag} Laboratoire MODAL'X, UMR CNRS 9023, Universit\'{e} Paris Nanterre, 200 avenue de la R\'{e}\-pu\-bli\-que, 92001 Nanterre, France.}
\email{yann.demichel@parisnanterre.fr}

\begin{document}

\begin{abstract}
A strange title, might you say: Answer is in the question! However, contrary to popular belief and numerous citations in the literature, the image of the snowflake curve is not present or even mentioned in von Koch's original articles. So, where and when did the first snowflake fall? Unravel the mystery of the snowflake curve with us on a journey through time.
\end{abstract}


\maketitle

\fussy

\section*{Introduction}\label{sec:intro}

\noi In 2024 we will commemorate the centenary of the death of Swedish mathematician Helge von Koch. We will also be celebrating the 120th anniversary of the birth of its famous curve. Less well known, von Koch was part of the first generation of mathematicians around Mittag-Leffler at the newly created {\it Stockholms Högskola} (see \cite{gar98}). While most of his work is now outdated, his name has gone down in history for having been given to one of the most emblematic geometric figures of the 20th century: the snowflake curve. Known for its a priori paradoxical property of possessing an infinite perimeter but delimiting a finite area, it owes its fame to its suggestive and particularly aesthetic shape. Yet a mystery remains. The snowflake curve is clearly based on von Koch's construction but it does not appear in either the original 1904 article or the extended 1906 version (see \cite{koc04,koc06}). In fact, it does not appear in any of the Swedish mathematician's works. So, who is behind the snowflake curve? When and where did it first appear? That is what we will try to find out.

\section{Follow tracks in the snow}\label{sec:Snowflake1}

\begin{figure}[h!]
\begin{center}
\includegraphics[width=0.55\textwidth]{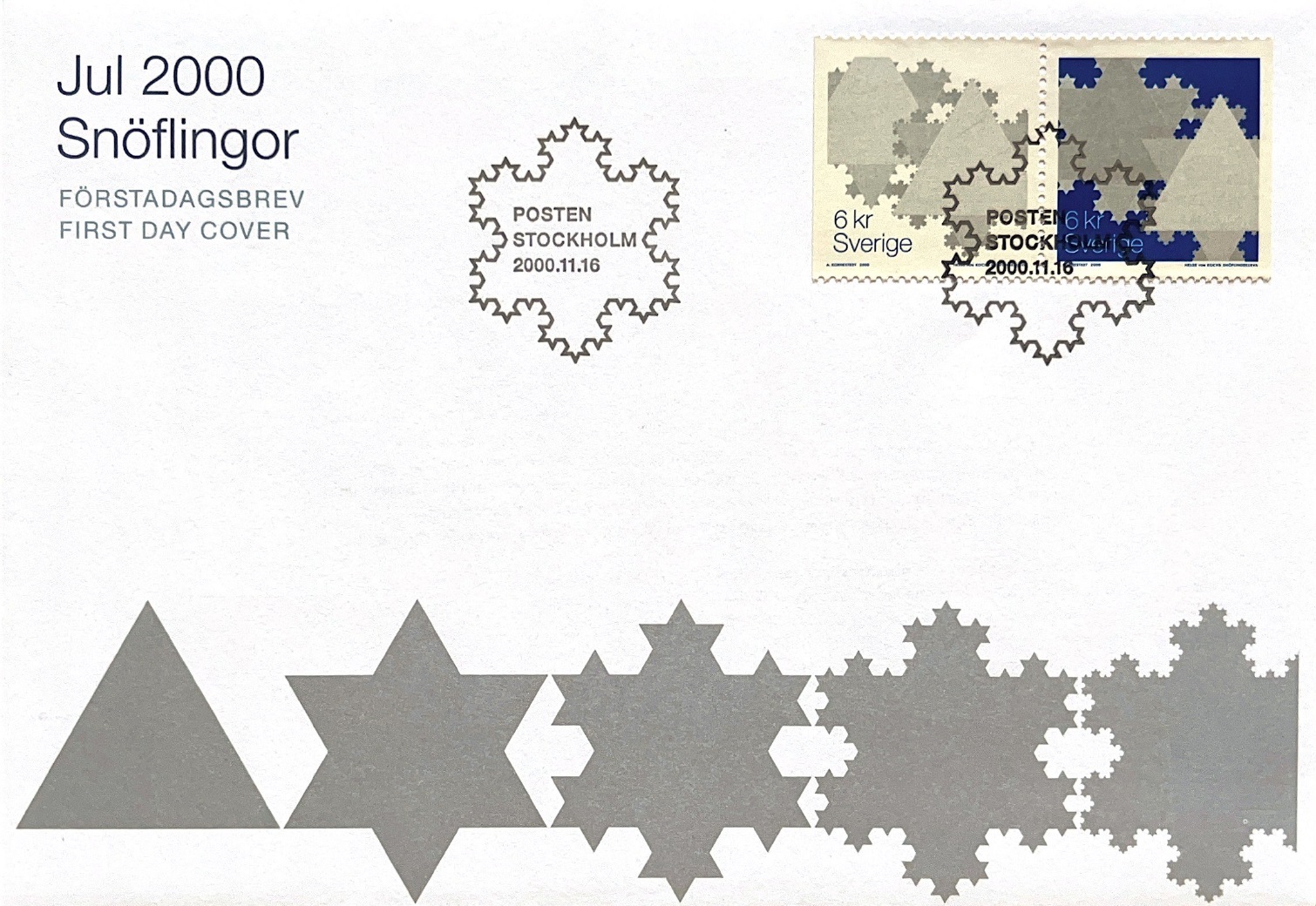}
\end{center}
\captionsetup{width=0.65\textwidth}
\caption{The limited-edition First Day Cover envelope printed by the Swedish Post Office on November 16, 2000, featuring the snowflake curve.}
\label{fig:Snowflake2000}
\end{figure}

\noi The starting point of our investigation is von Koch's native Sweden. Chosen as the emblem of the Swedish Mathematical Society, the snowflake curve is the pride of the profession. In 2000, declared the World Mathematical Year, the Swedish Post Office honored it for Christmas celebration. Two stamps featuring the snowflake are issued, affixed to an envelope showing the first stages of the snowflake construction, and canceled with a special postmark (see Figure \ref{fig:Snowflake2000}). The stamps are marked {\it Helge von Kochs snöflingekurva} and the back of the envelope reads:

\begin{quote}
``In 1904 the Swedish mathematician Helge von Koch (1870--1924) constructed an early example of fractals, the snowflake curve.''
\end{quote}

\noi Not surprisingly Sweden fully claims the property of the figure and its inventor, and invites us to follow the path of fractal geometry. A few years earlier it was in this context that the snowflake was twice the star of the science pages of the national daily newspaper {\it Svenska Dagbladet} (see \cite{bok83,bok91}). In 1991 it shared the bill with the original curve, demonstrating its obvious kinship. Two versions of the same curve, both attributed to von Koch. In the 1983 article, it is presented as the archetypal image of the new geometry developed by Benoît Mandelbrot.

\noi Mandelbrot said that von Koch's curve was an inspiration to him throughout his life. In his masterpiece {\it The Fractal Geometry of Nature}, he devotes an entire chapter to it and cites it no less than 200 times (see \cite{man82}). Because of their different fractal properties he was the first to clearly distinguish between the original curve, the snowflake curve and the filled snowflake. Three different fractal sets, all by one man: von Koch. Mandelbrot did not invent the snowflake. If that were the case he would have made it known. Moreover, he finds the name inappropriate, preferring `Koch Island', largely because he originally used von Koch's curve as a model for coastlines. This gives us an idea: What if the snowflake curve was born under a different name?

\ppass An article published in 1979 in {\it Le Petit Archimède}, a French magazine for young people devoted to mathematical recreation, gives us such an indication (see \cite{lpa79}). The snowflake curve is also called {\it van der Waerden's Thistle} because of its resemblance to the spiky flower, and may be the work of the Dutch mathematician. We have managed to trace this back to 1968. In his book {\it Stories About Sets}, the Russian mathematician Naum Vilenkin reproduced the snowflake curve and also attributed it to van der Waerden (see \cite[pp. 101--102]{vil68}). But in the same year in the Netherlands, the curve was presented in the mathematical journal for high school students {\it Pythagoras} only under the name {\it sneeuwvlok-kurve van Von Koch} (see \cite{pyt68}). We are on the wrong track here. The confusion probably stems from the proximity of the two family names due to the noble particle, and the fact that, in 1930, van der Waerden had exhibited a continuous but nowhere differentiable function whose construction process is somewhat reminiscent of von Koch's and whose graph does, indeed, closely resemble a thistle.

\ppass We decide to leave Europe and continue our investigation in the United States going back before Mandelbrot's work. Because Mandelbrot's work is authoritative, and because the construction is essentially identical to that of the 1904 curve, the snowflake has been definitively attributed to von Koch. If we want to know more about its inventor, we have to go back before the era of fractal geometry and look outside the realm of theoretical research.

\ppass Let us look at popularizing mathematics first. Martin Gardner, through his famous columns for the magazine {\it Scientific American}, has often described Mandelbrot's fractal objects, especially the snowflake curve. He always attributed it to von Koch, except in his first article in 1965, where the name is not mentioned (see \cite{gar65}). Instead, he amusedly pointed out that the snowflake curve was at the heart of a science-fiction novel written in 1956 by the British physiologist William Grey Walter. Might it be a literary invention? The title of the novel, {\it The Curve of the Snowflake}, could not be more eloquent. One of the characters tells his friend (see \cite[pp. 72--73]{wal56}):

\begin{quote}
``The boundary of a finite area can be a line of infinite length. [...] I can show you in a minute on a piece of notepaper. The snowflake curve, they call it. One of a family of pathological curves. Look!''
\end{quote}
and goes on to detail how to build the figure. Who is designated as {\it they}? The author does not say. Later in the novel the heroes are confronted with a mysterious vessel shaped like a three-dimensional snowflake. Thanks to its infinite surface area it can travel through time. This is the machine we need! The novel and its strange spaceship fascinated many readers. Joel Schneider was one of them.
He has often told the story of how he was puzzled by the strange machine he tried to draw when he was in junior high school. This led to his desire to become a mathematician. He made the generalization of the von Koch curve the subject of his first article, accepted in a mathematics journal when he was only 22 years old (see \cite{sch65}). In his article, the snowflake curve is simply called von Koch's curve, and the word snowflake is not used. Schneider cites von Koch's second article and only one other reference, a popular science book entitled {\it Mathematics and Imagination}, published in 1940.

\section{Let it snow, let it snow}\label{sec:Snowflake2}

\noi Co-authored by Edward Kasner, a professor at Columbia University, and James Newman, one of his students, the book is a potpourri of mathematical topics depicted and explained in an imaginative style. At the end of the last chapter on differential and integral calculus is an appendix entitled {\it Pathological Curves}. Precisely the expression used in Grey Walter's novel. The first patient examined by Kasner and Newman is nothing but the snowflake curve (see \cite[pp. 344--355]{kas40}). This gives its title to a subsection in which the authors give an explicit description of the iterative construction of the curve, starting with the equilateral triangle and illustrating it with the figure obtained after each of the first five stages. Then they ask the reader two questions: ``Why is it called the snowflake curve and why is it called pathological?'' Answers they gave are simple. The curve is named for its suggestive shape and its pathology is that it has an infinite perimeter although it can be drawn on a sheet of paper the size of a postage stamp. It is also mentioned that the curve has no tangent at any of its points. A footnote at the end of the chapter points out that the curve essentially provides a function that is continuous everywhere but differentiable nowhere. This clearly indicates that the authors were aware of von Koch's articles whose aim was precisely to give a simple and intuitive example of such a function. However they do not mention von Koch and his work nor do they attribute to him the principle of construction. Authors continue with three examples: first another variant of von Koch's curve they called {\it Anti-Snowflake Curve}, next a curve filling an entire cubical box, finally a curve called {\it Crisscross Curve} because of crossing itself at every one of its points. These two last curves copied Hilbert's and Sierpi\'{n}ski's respectively but their names are not mentioned. These deliberate omissions can be explained by the desire to maintain a certain mystery around these curves and to appeal to the general reader by preferring poetic names to the names of mathematicians. But it is also very tempting to see this as a way for the authors of claiming they invented both these variant shapes and their names.

\ppass This idea is attractive for at least two reasons. First, the book by Kasner and Newman made the snowflake widely known to teachers and the general public. Indeed, with the craze for mathematical popularization that emerged in those years, the book was a huge success in the United States, then in South America and Europe, including Sweden, where it was translated and widely distributed. Secondly, the book contains other terms that Kasner certainly invented, such as the word {\it googol} for the number formed by the digit 1 followed by 100 zeros.

\ppass We decide to follow this track and investigate Columbia's teaching faculty. In 1951, Howard Fehr, head of the Mathematics Education Department at {\it Teachers College}, published a training manual for future teachers that explained how to teach the mathematical concepts included in the secondary school mathematics curricula.
At the end of the chapter devoted to sequences, series and limits, Fehr argues that it is more relevant to give students concrete applications of exact computations rather than purely algebraic examples, and highlights the snowflake curve as a fun and valuable application (see \cite[pp. 140--142]{feh51}). Fehr indicates that this example is taken from various sources but writes that

\begin{quote}
``This particular type of limit was given its name by Professor Kasner of Columbia University, but it was originally suggested by the work of Helge von Koch and presented by Ludwig Boltzmann of Vienna in {\it Mathematische Annalen}, Volume 50, 1898, and called by him the H-curve.''
\end{quote}

\noi At last, we have a proof that the snowflake curve was named by Kasner. But the invention of the shape is in doubt. Was Boltzmann who constructed the first snowflake-like curve before von Koch and Kasner? The reference to the Austrian physicist is somewhat strange. In his article, Boltzmann gives a model for that he called the H-curve, a continuous but non-differentiable function, but its construction is not based on von Koch's geometric process. Moreover, von Koch makes no mention of this reference in his articles. This is clearly a wrong track again.

\section{The first snow}\label{sec:Snowflake3}

\noi Two questions remain now that we have identified the inventor of the snowflake curve: When did Kasner make this discovery and how did he arrive at it? The story of the invention of the {\it googol} provides valuable clues. We know that the term was not coined by Kasner himself in 1940 but by one of his young nephews during a walk in 1920. The discovery of the snowflake may therefore be older, and the name may have been suggested by a child in the family or a particularly imaginative pupil of Kasner.


\ppass Kasner completed his graduate studies at Columbia where he was one of the very first students to receive a doctorate in mathematics in 1899. After a postdoctoral year spent with Klein and Hilbert in Göttingen, he returned to the United States to embark on his teaching career at Barnard College, a liberal arts college affiliated with Columbia. Kasner's pedagogical approach was rather unconventional. His main aim was to develop a taste for mathematics in his students by stimulating their curiosity. To this end he never hesitated to surprise his audience with amusing and puzzling questions taken from natural observations: How many grains of sand make up the beach at Coney Island? How many drops of water fall on New York City on a rainy day? How long is the east coast of the United States? Kasner wanted to elicit spontaneous and original answers. The ones he considered the most outstanding were displayed on an `honor board' next to his office. This may be where the snowflake curve first appeared, as Jesse Douglas, another renowned student of Kasner's, testified (see \cite{dou58}):

\begin{quote}
``On the mathematics bulletin board near his office would be posted samples of his students' work which he considered of special merit or interest. Here might appear, for example, a carefully drawn figure of the `snowflake curve' (a certain continuous curve without a tangent) approximated by a patiently constructed polygon of 3,072 tiny sides.''
\end{quote}

\noi There is no doubt that this is our famous snowflake, shown in the fifth stage of its construction with its number of sides exactly equal to $3\times4^5$. The anecdote is undated but Douglas studied under Kasner from 1916 to 1920 and then taught at Columbia College from 1920 to 1926. This suggests a ten-year period for the discovery of the snowflake curve.

\begin{figure}[h!]
\begin{center}
\includegraphics[width=0.6\textwidth]{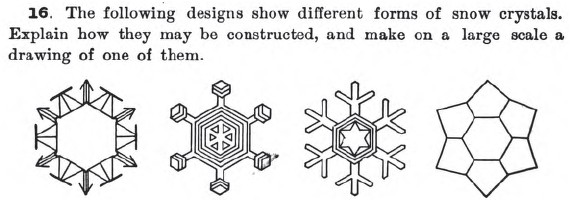}
\end{center}
\captionsetup{width=0.75\textwidth}
\caption{Drawings of snow crystals inspired by Bentley photographs. They illustrate an exercise in a famous geometry textbook used in American secondary schools around 1916.}
\label{fig:SnowCrystals1916}
\end{figure}

\ppass Kasner loved nature from which he drew his inspiration. Perhaps the idea of drawing a snowflake simply came to him during a walk in the snow. There is another less romantic but much more realistic hypothesis. Snowflakes were very fashionable in those years, popularized by the incredible photographs of a certain Wilson Bentley. Bentley was a farmer who lived in Jericho, Vermont, about 250 miles north of New York. Nicknamed `Snowflake', he had been fascinated by clouds, raindrops and snowflakes since childhood (see \cite{bla98}). In January 1885, at the age of 20, he became the first person to obtain a photomicrograph of a snow crystal. A self-taught researcher, Bentley has since accumulated hundreds of photographs which he has presented since 1902 in various journals and at public lectures, notably in New York. The remarkable structures revealed by its images have made snowflakes ideal objects to modernize the teaching of geometry. Rapidly, the images produced by Bentley or stylized representations of snowflakes appeared in American geometry textbooks. For example, the geometric construction exercise that Kasner gave his pupils was included in the famous 1916 book by Stone and Millis (see \cite[p. 252]{sm16} from which Figure \ref{fig:SnowCrystals1916} is taken).

\ppass These exercises undoubtedly inspired Kasner. By asking his students to draw a more complicated figure than the conventional ones he could discuss his favorite topic: the manifestation of very large numbers in nature. By starting with an equilateral triangle instead of a segment, von Koch's procedure made it possible to construct polygons with as many sides as desired still going with `natural' shapes.

\ppass Kasner must have been familiar with von Koch's curve. Indeed, he was fascinated by all these so-called crinkly curves. In a lecture on geometry given on September 24, 1904, at the St. Louis Congress of Arts and Sciences, he devoted a part to such curves, citing the historical examples of Weierstrass and Peano.
But he also talked about the curve that Osgood had constructed just the year before. This proves that Kasner was following this topic very closely. Thanks to Mittag-Leffler, von Koch's articles quickly circulated in the mathematical community. Although he could not quote von Koch in his lecture \--- his paper was published the month after the congress \--- it is clear that Kasner subsequently discovered this new curve with great interest.

\ppass It is therefore possible that Kasner's snowflake could have been designed and used before 1916. Several pieces of evidence are in favour of this. First Douglas may not have seen the snowflake drawing pinned on the famous board with his own eyes. Someone might have told him a story that happened a few years ago. Secondly, if the geometry exercise did not appear in the first edition of the Stone-Millis book in 1910 \--- it was added to the revised edition of 1916 \---, a similar exercise can be found in another reference book published in 1912 (see \cite[p. 55]{bet12}).

\ppass This new date allows us to consider a final hypothesis, namely that Kasner and von Koch met and exchanged views on the snowflake curve. Indeed, the two men were both invited to address the 1912 International Congress of Mathematicians in Cambridge, England.
Although eight years younger than von Koch, Kasner has a lot in common with him. Professionally, he was a brilliant and precocious mathematician, influenced by the work and ideas of Poincaré, who emphasized the value and importance of intuition, particularly in teaching. On a more personal level, Kasner and von Koch shared a passion for art, a love of nature and an affection for children. This had a major influence on their approach to mathematics and teaching. At the time, von Koch was Sweden's delegate to the International Mathematical Instruction Commission. Reports from all the national commissions were to be presented in the didactic section of the congress. Even if the Swedish commission's reports were not read by von Koch \--- he attended the analysis session held at the same time \---, Kasner could not have been unaware of his colleague's interest in teaching. Everything indicates that Kasner might have wanted to talk to von Koch. Could they have discussed the snowflake curve? Could it have been von Koch who suggested this variant to Kasner? We can only speculate here. At the end of our investigation, the snowflake curve finally did not reveal all its secrets...

\bibliographystyle{plain}

\end{document}